\documentclass{article}
\usepackage{latexsym}
\usepackage{amssymb}
\usepackage{amsmath}
\usepackage{layout}
\newtheorem{teo}{Theorem}[section]
\newtheorem{lema}[teo]{Lemma}

\newtheorem{prop}[teo]{Proposition}
\newtheorem{obs}[teo]{Remark}

\newenvironment{dem}{\text\bf Proof:}{}
\newenvironment{dem2}{}{}
\newcommand{\R}{\mathbb{R}}
\newcommand{\1}{{\bf 1}}

\newcommand{\Pro}{\mathbb{P}}
\newcommand{\cL}{{\cal L}}

\newcommand{\al}{\alpha}
\newcommand{\la}{\lambda}
\newcommand{\si}{\sigma}
\begin{document}

\begin{center}
\huge {

Convergence of delay differential equations
driven by fractional Brownian motion
}

\vspace{.5cm}

\normalsize {\bf Marco Ferrante$^{1}$} and {\bf Carles
Rovira$^{2,*}$}

{\footnotesize \it $^1$ Dipartimento di Matematica Pura ed App.,
Universit\`a di Padova,
Via Trieste 63, 35121-Padova, Italy.

$^2$ Facultat de Matem\`atiques, Universitat de Barcelona, Gran
Via 585, 08007-Barcelona, Spain.

 {\it E-mail addresses}: ferrante@math.unipd.it,
Carles.Rovira@ub.edu}

{$^{*}$corresponding author}

\end{center}

\begin{abstract}%
In this note we prove an existence and uniqueness result of
solution for stochastic differential delay equations with
hereditary drift driven by a fractional Brownian motion with Hurst
parameter $H > 1/2$. Then, we show that, when the delay goes to
zero, the solutions to these equations converge, almost surely and
in $L^p$, to the solution for the equation without delay. The
stochastic integral with respect to the fractional Brownian motion
is a pathwise Riemann-Stieltjes integral.

\end{abstract}



\setcounter{section}{-1}

\section {Introduction}

Consider the stochastic
differential equation on $\R^d$
\begin{eqnarray}\label{equau}
X^r (t) &= &\eta (0) + \int_0^t b(s, X^r) ds + \int_0^t \si(s, X^r(s-r)) dW_s^H, \quad t \in (0,T],\cr
X^r (t) &= &\eta (t), \qquad t \in [-r,0].
\end{eqnarray}
Here $r$ denotes a strictly positive time delay,
$W^H=\{W^{H,j},j=1,\ldots,m\}$ are independent fractional
Brownian motions with Hurst parameter $H > \frac12$ defined
in a complete probability space $(\Omega, {\cal F}, \Pro)$,
$b(s,X^r)$, the hereditary term, depends on the path
$\{X^r(u), -r \le u \le s \}$, while $\eta: [-r, 0] \to \R^d$ is a
smooth function.
We call (\ref{equau}) a delay differential equations
with hereditary drift driven by a fractional Brownian motion and
to the best of our knowledge this problem has not been
considered before in the wide
literature on stochastic differential equations.

As usual in this field, we have to specify how we intend
the stochastic integral in (\ref{equau}),
being its definition not unique.
Since $H > \frac12$, we can define the integral with respect to
fractional Brownian motion using a pathwise approach.
Indeed, if we have a stochastic processes $\{u(t), t \ge 0\}$
whose trajectories are $\lambda$-H\"older continuous with $\lambda
> 1-H$, then the Riemann-Stieltjes integral $\int_0^T u(s) dW_s^H$
exists for each trajectory (see Young \cite{Y}). Using the
techniques introduced by Young \cite{Y} and the $p$-variation
norm, Lyons \cite{L} began the study of integral equations driven
by functions with bounded $p$-variation, with $p\in [1,2)$.
Then Zahle \cite{Z} introduced a generalized Stieltjes integral using
the techniques of fractional calculus. The integral is expressed
in terms of fractional derivative operators and it coincides with
the Riemann-Stieltjes integral $\int_0^T f dg$ when the functions
$f$ and $g$ are H\"older continuous of orders $\lambda$ and $\beta$,
respectively, with $\lambda + \beta > 1$. Using this
Riemann-Stieltjes integral, Nualart and Rascanu \cite{NR} obtained
the existence and uniqueness of solution for a class of integral
equations without delay and they also proved that the solution is
bounded on a finite interval.

In our paper, using also the Riemann-Stieltjes integral, we will
first prove the existence and uniqueness of a solution to
equation (\ref{equau}),
extending the results in Nualart and Rascanu \cite{NR}.
Then we will study the convergence of the solutions
of these equations when the delay $r$ tends to zero
and the drift coefficient $b$ depends on $(s,X^r(s))$.
As occurs for the Brownian motion case, we are able to prove that
the solution to the delay equation converges,
almost surely and in $L^p$, to the
solution of the equation without the delay.
All along the paper,
we will prove first our results for
deterministic equations and then we will easily
apply them pathwise to fractional Brownian motion.

There are many references on stochastic systems with delay (see
for instance \cite{M}), but the literature about stochastic
differential equations with delay driven by a fractional Brownian
motion is scarce. In a previous paper
\cite{FR} we obtain the
existence and uniqueness of solution and the smoothness of the
density when $H>1/2$ under strong hypothesis, using only
techniques of the classical stochastic calculus.
That approach is unfortunately not suitable for further investigation,
like the presence of an hereditary drift and the convergence
when the delay tends to zero.
Using rough path analysis,
Neuekirch, Nourdin and Tindel \cite{NNT} considered the case $H >
1/3$. Nowadays, Le\'on and Tindel \cite{LT} are studying the
existence of solution and its regularity when $H> 1/2$.

The structure of the paper is as follows:
in the next section we state the main results of our paper.
In Section 2 we give some useful estimates for Lebesgue
integrals and for Riemann-Stieltjes ones.
Section 3 is devoted to obtain the existence,
uniqueness and boundedness for the solution for deterministic equations.
Section 4 contains the study of the convergence of the deterministic equations.
In Section 5 we apply the results of the previous sections to
stochastic equations driven by fractional Brownian motion and
we give the proofs of our main theorems.
Finally, in Section 6 we recall a couple of technical results.

We will denote by $C_\alpha$ a constant that will change from line
to line.

\section {Main results}

Let $\al \in (\frac12,1)$ and $r>0$. We will denote by $W_0^{\al,\infty} (-r,T;\R^d)$ the space of mesurable functions $f: [0,T] \to \R^d$ such that
\begin{eqnarray*}
\Vert f \Vert_{\alpha,\infty(r)} := \sup_{t \in [-r,T]} \Big( \vert
f(t)\vert + \int_{-r}^t \frac{\vert f(t)-f(s) \vert}{
(t-s)^{\alpha+1}}ds \Big)< \infty.
\end{eqnarray*}
For any $\lambda \in (0,1]$, we will consider $C^\lambda (-r,T; \R^d)$ the space of $\lambda-$H\"older continuous functions $f: [0,T] \to \R^d$ such that
\begin{eqnarray*}
\Vert f \Vert_{\lambda(r)} := \Vert f \Vert_{\infty(r)} + \sup_{-r
\le s < t \le T}  \frac{\vert f(t)-f(s) \vert}{ (t-s)^{\lambda}}<\infty,
\end{eqnarray*}
where
\begin{eqnarray*}
\Vert f \Vert_{\infty(r)} := \sup_{s \in [-r,T]} \vert f(s) \vert.
\end{eqnarray*}
Note that when $r=0$, we shall omit $(r)$ in the name of the
corresponding norm.

We also need to consider the spaces $W_0^{\al,\infty}(-r,0;\R^d)$ and $C^{\la} (-r,0; \R^d)$
with the corresponding norms
\begin{eqnarray*}
&& \Vert \eta \Vert_{\al,\infty(-r,0)}:=  \sup_{t \in [-r,0]}
\Big( \vert \eta(t)\vert + \int_{-r}^t \frac{\vert \eta(t)-\eta(s)
\vert}{ (t-s)^{\alpha+1}}ds \Big),
 \cr
 &&\Vert \eta \Vert_{\lambda(-r,0)} := \Vert \eta
\Vert_{\infty(-r,0)} + \sup_{-r \le s < t \le 0}  \frac{\vert
\eta(t)-\eta(s) \vert}{ (t-s)^{\lambda}} \qquad {\rm where} \quad
\Vert \eta \Vert_{\infty(-r,0)} := \sup_{s \in [-r,0]} \vert
\eta(s) \vert.
\end{eqnarray*}

\bigskip

Let us consider the following hypothesis:

\begin{itemize}

\item{\bf (H1)} $\si: [0,T] \times \R^d \to \R^d \times \R^m$ is a
measurable function such that  $\si(t,x)$ is differentiable in $x$ and there
exists some constants $0 < \beta, \delta \le 1$ and for every $N
\ge 0$ there exists $M_N >0$ such that the following properties
hold:
\begin{enumerate}
\item
$ \vert \si(t,x)-\si(t,y) \vert \le M_0 \vert x-y \vert, \quad
\forall x,y \in \R^d, \quad \forall t \in [0,T]$,

\item
$ \vert \partial_{x_{i}} \si(t,x) -  \partial_{y_{i}}
\si(t,y)\vert \le M_N \vert x-y \vert^\delta, \quad \forall \vert
x \vert , \vert y \vert \le N, \quad \forall t \in [0,T]$, for
each $i=1,...,d.$

\item
$ \vert \si(t,x)-\si(s,x) \vert +   \vert \partial_{x_{i}}
\si(t,x) -  \partial_{x_{i}} \si(s,x)\vert             \le M_0
\vert t-s \vert^\beta, \quad \forall x\in \R^d, \quad \forall t, s
\in [0,T]$ for each $i=1,...,d.$

\end{enumerate}

\item {\bf (H2)} $b: [0,T] \times C(-r,T;\R^d) \to \R^d$ is a
measurable function such that for every $t >0$ and $f \in C(-r,T;\R^d)$,
$b(t,f)$ depends only on $\{f(s); -r \le s \le t\}$. Moreover, there exists $b_0 \in L^\rho(0,T;\R^d)$ with $\rho
\ge 2$ and $\forall N \ge 0$ there exists $L_N>0$ such that
\begin{enumerate}
\item
$ \vert b(t,x)-b(t,y) \vert \le L_N \sup_{-r \le s \le t} \vert
x(s)-y(s) \vert, \quad  \forall x,y ,\, \Vert x \Vert_{\infty(r)}
\le N, \Vert y \Vert_{\infty(r)} \le N, \forall t \in [0,T]$,

\item
$ \vert b(t,x)\vert \le L_0 \sup_{-r \le s \le t} \vert x(s) \vert +
b_0(t),\quad \forall t \in [0,T]$
\end{enumerate}

\item {\bf (H3)} There exists $\gamma \in [0,1]$ and $K_0 >0$ such
that
$$
\vert \si(t,x) \vert \le K_0 ( 1 + \vert x \vert^\gamma), \quad
\forall x \in \R^d, \forall t \in [0,T].$$

\end{itemize}

\bigskip

Under these assumption we are able to prove that our problem admits
a unique solution.
The result of existence and uniqueness reads as follows:

\medskip

\begin{teo}\label{teoexis}
Assume that $\eta \in W_0^{\al,\infty}(-r,0;\R^d) \cap C^{1-\al}
(-r,0; \R^d)$ and that $b$ and $\si$ satisfy hypothesis {\bf (H1)}
and  {\bf (H2)} with $\beta > 1 - H, \delta > \frac{1}{H} -1$. Set
$$ \al_0:=\min \{\frac12, \beta, \frac{\delta}{1+\delta} \}.$$
Then if $\al \in (1-H, \al_0)$  and $\rho \le \frac{1}{\al}$,
equation (\ref{equau}) has an unique solution $$X^r \in
L^0(\Omega, {\cal F}, \Pro; W_0^{\al,\infty} (-r,T; \R^d))$$  and
for $P$-almost all $\omega \in \Omega,$ $X^r(\omega,.) \in
C^{1-\al}(-r,T;\R^d)$.

Moreover, if $\alpha \in (1-H, \al_0 \vee (2-\gamma)/4)$ and {\bf
(H3)} holds then $E ( \Vert X \Vert^p_{\al,\infty(r)} ) < \infty
\forall p \ge 1.$
\end{teo}

\bigskip

In order to study the convergence when the delay goes to zero,
we will consider the particular case of our initial equation
(\ref{equau}), where the coefficient $b$
does not depend on the whole trajectory,

\begin{eqnarray}\label{equags}
X^r (t) &= &\eta (0) + \int_0^t b(s, X^r(s)) ds + \int_0^t \si(s,
X^r(s-r)) dW_s^H, \quad t \in (0,T],\cr X^r (t) &= &\eta (t),
\qquad t \in [-r,0].
\end{eqnarray}

We will assume that $b$ satisfies the new set of hypothesis:

{\bf (H2')} $b:[0,T] \times \R^d \to \R^d$ is a measurable
function such that there exists $b_0 \in L^\rho(0,T;\R^d)$ with
$\rho \ge 2$ and $\forall N \ge 0$ there exists $L_N>0$ such that
\begin{enumerate}
\item
$ \vert b(t,x)-b(t,y) \vert \le L_N  \vert x-y \vert, \quad
\forall x,y ,\, \vert x \vert \le N,  \vert y \vert \le N, \quad
\forall t \in [0,T],$

\item
$ \vert b(t,x)\vert \le L_0  \vert x(t) \vert + b_0(t), \quad
\forall t \in [0,T].$
\end{enumerate}

\medskip

Denoting by $X$ the solution of the
stochastic differential equation on $\R^d$ without delay:
\begin{equation*}
X (t) = \eta (0) + \int_0^t b(s, X(s)) ds + \int_0^t \si(s, X(s))
dW_s^H, \quad t \in [0,T].
\end{equation*}
we are able to prove the following result:

\smallskip

\begin{teo}\label{convstoch}
Assume that $\eta \in W_0^{\al,\infty}(-r_0,0;\R^d) \cap C^{1-\al}
(-r_0,0; \R^d)$ and that $b$ and $\si$ satisfy hypothesis {\bf
(H1)} and  {\bf (H2')} with $\beta > 1 - H, \delta > \frac{1}{H}
-1$. Set
$$ \al_0:=\min \{\frac12, \beta, \frac{\delta}{1+\delta} \}.$$
Then if $\al \in (1-H, \al_0)$  and $\rho \le \frac{1}{\al}$ for
$P$-almost all $\omega \in \Omega$
$$
\lim_{r \to 0} \Vert X(\omega,.) - X^r(\omega,.) \Vert_{\al,\infty} = 0.
$$
Moreover, if $\alpha \in (1-H, \al_0 \vee (2-\gamma)/4)$ and {\bf
(H3)} holds then
$$
\lim_{r \to 0} E( \Vert X - X^r \Vert_{\al,\infty}^p ) = 0, \quad \forall p \ge 1.
$$

\end{teo}

\section{Estimates for the integrals}

In this section we will obtain some estimates for the Lebesgue
integral and for the pathwise Riemann-Stieltjes integral. In both
cases, we will recall some well-known results and we will obtain
some estimates well posed to our equations.
\medskip
In the space $W_0^{\al,\infty} (0,T;\R^d)$ we need to introduce a new norm, that is, for any $\lambda \ge 1$
\begin{eqnarray*}
\Vert f \Vert_{\alpha,\lambda(r)} := \sup_{t \in [-r,T]} \exp( -
\lambda t) \Big( \vert f(t) \vert + \int_{-r}^t \frac{\vert
f(t)-f(s) \vert}{ (t-s)^{\alpha+1}}ds \Big).
\end{eqnarray*}
It is easy to check that, for any $\lambda \ge 1$, this norm is
equivalent to $\Vert f \Vert_{\alpha,\infty(r)}.$

\subsection{The Lebesgue integral}

Let us consider first the ordinary Lebesgue integral. Given
$f:[0,T] \to \R^d$ a measurable function we define
$$ F(f) (t) = \int_0^t f(s)  ds.$$

We recall first a result from \cite{NR} (see Proposition 4.3).

\begin{prop}\label{Prop43nua}
Let $0 < \alpha <  \frac12$ and $f:[0,T] \to \R^d$ be a measurable
function. If $$\sup_{t \in [0,T]} \int_0^t \frac{\vert f(s)
\vert}{(t-s)^\al}<\infty$$ then $F(f)(.) \in
W_0^{\al,\infty}(0,T;\R^d)$ and
\begin{equation*}
\vert F(f)(t)\vert + \int_0^t \frac{\vert F(f)(t)-F(f)(s) \vert}{(t-s)^{\al+1}} ds \le C_{\al,T} \int_0^t
\frac{\vert f(s) \vert}{(t-s)^\al} ds.
\end{equation*}
\end{prop}

\bigskip

Given
$f:[0,T] \to \R^d$ a measurable function we define
$$ F^{(b)}(f) (t) = \int_0^t b(s,f)  ds.$$

\begin{prop}\label{desleb}
Assume that $b$ satisfies {\bf (H2)} with $\rho=\frac{1}{\al}$. If $f \in W_0^{\al,\infty}(-r,T;\R^d)$ then
$F^{(b)} (f) (.) = \int_0^. b(s,f) ds \in C^{1-\al} (0,T;\R^d)$ and
\begin{eqnarray*}
&& i) \ \ \Vert F^{(b)}(f)\Vert_{1-\al}  \le d^{(1)} (1+  \Vert f
\Vert_{\infty(r)}), \cr && ii) \ \ \Vert F^{(b)}(f)\Vert_{\al,\la}
\le d^{(2)} \left(\frac{1}{\la^{1-2\al}}+  \frac{\Vert f
\Vert_{\al,\la(r)}}{\la^{1-\al}} \right) \le
\frac{d^{(2)}}{\la^{1-2\al}} \left(1+ \Vert f \Vert_{\al,\la(r)}
\right),
\end{eqnarray*}
for all $\la \ge 1$ where $d^{(i)}, i \in \{1,2\}$ are positive constants depending only on $\al, T, L_0$
and $B_{0,\al} = \Vert b_0 \Vert_{L^{1/\al}}$.

If $f, h  \in W_0^{\al,\infty}(-r,T;\R^d)$ such that $\Vert f \Vert_{\infty(r)} \le N, \Vert h \Vert_{\infty(r)} \le N,$ then
\begin{equation*}
\Vert F^{(b)}(f) - F^{(b)}(h) \Vert_{\al,\la}  \le \frac{d_N}{\la^{1-\al}} \Vert f-h \Vert_{\al,\la(r)}
\end{equation*}
for all $\la \ge 1$ where $d_N=C_{\al,T} L_N \Gamma(1-\al)$ depends on $\al, T$ and $L_N$ from {\bf (H2)}.
\end{prop}

\begin{dem}
It follows the ideas of Proposition 4.4 in \cite{NR}. For the sake
of completeness, we will give a sketch of the proof.

In order to simplify the presentation, we will assume $d=1$.
For $f \in W_0^{\al,\infty}(-r,T)$ and $0 \le s \le t \le T$ we can write
\goodbreak
\begin{eqnarray*}
\vert F^{(b)}(f)(t) - F^{(b)}(f)(s) \vert & \le & \int_s^t ( L_0
\sup_{-r \le v \le u} \vert f(v) \vert + b_0(u)) du \cr & \le &
L_0 \Vert f \Vert_{\infty(r)} (t-s) + \int_s^t b_0(u)du \cr & \le
& ( L_0 T^\al \Vert f \Vert_{\infty(r)} + B_{0,\al}
)(t-s)^{1-\al},
\end{eqnarray*}
where $B_{0,\al} := \Vert b_0 \Vert_{L^{1/\al}}$. So $ F^{(b)}(f) \in C^{1-\al} (0,T)$ and $i)$
is true with $$d^{(1)}=(L_0 T^\al + B_{0,\al})(1+T^{1-\al}).$$

On the other hand, from Proposition \ref{Prop43nua},
\begin{eqnarray*}
&& \vert F^{(b)}(f)(t) \vert + \int_0^t \frac{\vert F^{(b)}(f)(t) - F^{(b)}(f)(s) \vert}{(t-s)^{\al+1}}
 ds \cr
 &  & \qquad \le C_{\al,T}  \int_0^t \frac{L_0 \sup_{-r \le u \le s} \vert f(u) \vert + b_0(s)}{(t-s)^{\al}}
 ds \cr
&  & \qquad \le C_{\al,T}  \left( L_0 \int_0^t \frac{ \sup_{-r \le u \le s} \vert f(u) \vert }{(t-s)^{\al}}
 ds + \left(\int_0^t (t-s)^{-\al/(1-\al)} ds \right)^{1-\al}   B_{0,\al} \right).
\end{eqnarray*}
Using that
$$ e^{-\la s} \sup_{-r \le u \le s} \vert f(u) \vert \le  \sup_{-r \le u \le s} \vert f(u) \vert e^{-\la u},$$
we obtain that
\begin{eqnarray*}
&& \sup_{t \in [0,T]}  e^{-\la t}  \int_0^t \frac{\sup_{-r \le u \le s} \vert f(u) \vert}{(t-s)^{\al}}
 ds  \le \sup_{t \in [0,T]}   \int_0^t \frac{e^{-\la (t-s)}}{(t-s)^{\al}}  \sup_{-r \le u \le s} \vert f(u) e^{-\la u} \vert
 ds \cr
&  & \qquad \le \sup_{-r \le u \le T}  \vert f(u) \vert  e^{-\la
u}
 \sup_{t \in [0,T]}   \int_0^t \frac{e^{-\la (t-s)}}{(t-s)^{\al}}
 ds \le \la^{\al-1} \Gamma(1-\al) \sup_{-r \le u \le T}  \vert f(u) \vert e^{-\la u}.
\end{eqnarray*}
So, using that
$$\int_0^t \frac{e^{-\lambda(t-s)}}{(t-s)^\alpha} ds \le
\lambda^{\alpha-1} \Gamma(1-\alpha)\quad {\rm and} \quad \sup_{t
\in [0,T]} t^\mu e^{-\lambda t} \le \left( \frac{\mu}{\lambda}
\right)^\mu e^{-\mu},$$ for all $\la \ge 1$
\begin{eqnarray*}
&& \Vert F^{(b)}(f)\Vert_{\al,\la}  \le C_{\al,T} L_0 \la^{\al-1}
\Gamma(1-\al) \sup_{-r \le s \le T}  \vert f(s) e^{-\la s} \vert +
C_{\al,T}  \frac{(1-\al)^{1-\al}}{(1-2\al)^{\al}} e^{2\al-1}
B_{0,\al} \la^{2\al-1} \cr &&\qquad \le d^{(2)} ( \la^{2\al-1}  +
\la^{\al-1} \Vert f \Vert_{\al,\lambda(r)}),
\end{eqnarray*}
and {\it ii)} becomes true with
$$ d^{(2)} = C_{\al,T} [ L_0  \Gamma(1-\al) +
 (1-2\al)^{-\al} (1-\al)^{1-\al} e^{2\al-1} B_{0,\al}].
$$

Finally, if $f, h  \in W_0^{\al,\infty}(-r,T)$ such that $\Vert f \Vert_{\infty(r)} \le N, \Vert h \Vert_{\infty(r)} \le N,$ using similar computations we obtain that
\begin{eqnarray*}
\Vert F^{(b)}(f)-F^{(b)}(h) \Vert_{\al,\la} & \le & C_{\al,T}
\sup_{t \in [0,T]} e^{-\la t}  \int_0^t \frac{L_N  \sup_{-r \le u
\le s} \vert f(u)-h(u) \vert}{(t-s)^\al} ds \cr & \le & C_{\al,T}
L_N \sup_{-r \le u \le T} \left( e^{-\la u}  \vert f(u)-h(u) \vert
\right) \sup_{t \in [0,T]} \int_0^t \frac{e^{-\la
(t-s)}}{(t-s)^{\al}}
 ds
\cr
& \le & \frac{1}{\la^{1-\al}} d_N \Vert f - h \Vert_{\al,\la(r)},
\end{eqnarray*}
with $d_N=C_{\al,T} L_N \Gamma(1-\al).$
\hfill $\Box$
\end{dem}

%

\begin{obs}\label{obsb}
If we assume
{\bf (H2')} and define
$$ F^{(b)}(f) (t) = \int_0^t b(s,f(s))  ds,$$
we obtain a version of Proposition \ref{desleb} with $r=0$ that is
an extension of Proposition 4.4 in \cite{NR}.
\end{obs}

\subsection{The Riemann-Stieltjes integral}

The Riemann-Stieltjes integral introduced by Z\"ahle \cite{Z} is based on fractional integrals
and derivatives (see \cite{SKM}). We will refer the reader to the papers of Z\"ahle \cite{Z} and
Nualart and Rascanu \cite{NR} for the general  theory.

Here, we will just recall some basic results. Fixed a parameter $0
< \al < 1/2$, let us consider $W_T^{1-\al,\infty}(0,T)$ the space
of measurable functions $g:[0,T] \to \R$ such that
\begin{eqnarray*}
\Vert g \Vert_{1-\alpha,\infty,T} := \sup_{0<s<t<T} \Big( \frac{\vert
g(t)-g(s)\vert}{(t-s)^{1-\al}} + \int_{0}^t \frac{\vert g(u)-g(s) \vert}{
(u-s)^{2-\alpha}}du \Big)< \infty.
\end{eqnarray*}
If $g \in W_T^{1-\al,\infty}(0,T)$, we can define
$$
\Lambda_\al(g):= \frac{1}{\Gamma(1-\al)}
\sup_{0<s<t<T} |(D^{1-\al}_{t-} g_{t-})(s)|
$$
where $\Gamma(\cdot)$ is the Euler function and
$D^{1-\al}_{t-}$ denotes the Weyl derivative
(see Nualart and Rascanu \cite{NR} for more details).
We get
$$\Lambda_\al(g) \le \frac{1}{\Gamma(1-\al)\Gamma(\al)} \Vert g \Vert_{1-\alpha,\infty,T} < \infty.
$$
We also consider $W_0^{\al,1}(0,T)$ the space of measurable functions $f:[0,T] \to \R$ such that
\begin{eqnarray*}
\Vert f \Vert_{\alpha,1} := \int_0^T \frac{\vert
f(s)\vert}{s^{\al}} ds+ \int_0^T  \int_{0}^s \frac{\vert f(s)-f(u)
\vert}{ (s-u)^{\alpha+1}}du ds< \infty.
\end{eqnarray*}
Given two functions $g
\in W_T^{1-\al, \infty} (0,T)$ and $f \in W_0^{\al, 1} (0,T)$,
we can define
$$G(f) (t):= \int_0^t f(s) dg_s= \int_0^T f(s) \1_{(0,t)} (s) dg_s.$$
It holds that
$$\vert \int_0^t f(s) dg_s \vert \le   \Lambda_\al(g) \Vert f \Vert_{\al,1}.$$

The following estimates are proved in \cite{NR} (see Proposition
4.1)
\begin{equation}\label{NR418}
\vert G(f)(t) \vert\le \Lambda_\al(g)  \Big( \int_0^t \frac{\vert
f(s)\vert}{s^\al} ds + \alpha \int_{0}^t \int_0^s\frac{\vert
f(s)-f(y) \vert}{ (s-y)^{\alpha+1}}dyds \Big),
\end{equation}

\begin{equation}\label{NR417}
\int_0^t \frac{\vert G(f)(t) - G(f)(s) \vert}{ (t-s)^{\alpha+1}}
ds  \le  \Lambda_\al(g)  \Big( C_\al \int_0^t \frac{\vert
f(s)\vert}{(t-s)^{2\al}} ds +
 \int_{0}^t \int_0^s
\frac{\vert f(s)-f(y) \vert}{ (s-y)^{\alpha+1}} (t-y)^{- \al} dyds
\Big).
\end{equation}

Let us consider the term
$$
G_r^{(\sigma)} (f)(t) = \int_0^t \sigma (s, f(s-r) ) dg_s.
$$

\begin{prop}\label{desint}
Assume that $\si$ satisfies
{\bf (H1)}. If $f \in W_0^{\al,\infty}(-r,T;\R^d)$ then
$$
G^{(\si)}_r (f) \in C^{1-\al}(0,T;\R^d) \subset
W_0^{\al,\infty}(0,T;\R^d),
$$
and
\begin{eqnarray*}
&& i) \ \ \Vert G^{(\si)}_r (f) \Vert_{1-\al}  \le \Lambda_\al(g)
d^{(3)} (1+  \Vert f \Vert_{\al,\infty(r)}), \cr
&& ii) \ \ \Vert
G^{(\si)}_r (f) \Vert_{\al,\la}  \le \frac{\Lambda_\al(g)
d^{(4)}}{\la^{1-2\al}} (1+  \Vert f \Vert_{\al,\la(r)}),
\end{eqnarray*}
for all $\la \ge 1$ where $d^{(i)}, i \in \{3,4\}$ are positive constants independent of $\la, f$ and $g$.

If $f, h  \in W_0^{\al,\infty}(-r,T;\R^d)$ such that $\Vert f \Vert_{\infty(r)} \le N, \Vert h \Vert_{\infty(r)} \le N,$ then
\begin{equation*}
\Vert G^{(\si)}_r (f) - G^{(\si)}_r (h) \Vert_{\al,\la}  \le
\frac{\Lambda_\al(g) d^{(2)}_N}{\la^{1-2\al}}
(1 + \Delta_r(f)+\Delta_r(h)) \Vert f-h \Vert_{\al,\la(r)}
\end{equation*}
for all $\la \ge 1$ where
$$
\Delta_r(f)=\sup_{-r \le u \le T} \int_{-r}^u \frac{ \vert f(u)-f(s) \vert^\delta}{(u-s)^{\al+1}} ds,$$
and $d_N^{(2)}$ does not depend on $\la$ and $g$.
\end{prop}

\begin{dem}
This proposition is a consequence of Proposition 4.2 in \cite{NR}.
Given $f \in W_0^{\al,\infty}(-r,T;\R^d)$ let us define $f^*:[0,T]
\to \R^d$ such that $f^*(s):=f(s-r)$. Clearly $f^* \in
W_0^{\al,\infty}(0,T;\R^d)$ and
\begin{eqnarray*}
\Vert f^* \Vert_{\alpha,\infty} & = & \sup_{t \in [0,T]} \Big(
\vert f^*(t)\vert + \int_{0}^t \frac{\vert f^*(t)-f^*(s) \vert}{
(t-s)^{\alpha+1}}ds \Big)\cr & = & \sup_{t \in [0,T]} \Big( \vert
f(t-r)\vert + \int_{0}^t \frac{\vert f(t-r)-f(s-r) \vert}{
(t-s)^{\alpha+1}}ds \Big)\cr & = & \sup_{u \in [-r,T-r]} \Big(
\vert f(u)\vert + \int_{-r}^{u} \frac{\vert f(u)-f(s) \vert}{
(u-s)^{\alpha+1}}ds \Big) \le  \Vert f
\Vert_{\alpha,\infty(r)},
\end{eqnarray*}
and
\begin{eqnarray*}
\Vert f^* \Vert_{\alpha,\lambda} & = & \sup_{t \in [0,T]} e^{ -
\lambda t} \Big( \vert f^*(t) \vert + \int_{0}^t \frac{\vert
f^*(t)-f^*(s) \vert}{ (t-s)^{\alpha+1}}ds \Big)\cr
 & = & e^{-\la r} \sup_{u \in [-r,T-r]} e^{ -
\lambda u} \Big( \vert f(u) \vert + \int_{-r}^u \frac{\vert
f(u)-f(s) \vert}{ (u-s)^{\alpha+1}}ds \Big)\cr
& \le & e^{-\la r}
\Vert f \Vert_{\alpha,\lambda(r)}
\end{eqnarray*}
and finally
$$
\Delta_0(f^*)=\sup_{0 \le u \le T} \int_{0}^u \frac{ \vert f^*(u)-f^*(s) \vert^\delta}{(u-s)^{\al+1}} ds
=\sup_{-r \le v \le T-r} \int_{-r}^v \frac{ \vert f(v)-f(s) \vert^\delta}{(v-s)^{\al+1}} ds
\le \Delta_r(f).
$$
Then, we only have to apply Proposition 4.2 in \cite{NR} to $f^*$.

\hfill $\Box$
\end{dem}

\begin{obs}\label{obssi}
If  we consider the case when $r=0$, we obtain a version of
Proposition \ref{desint} that coincides with Proposition 4.2 in
\cite{NR}.
\end{obs}

We need an additional estimate in order to be able
to bound the norms of the solutions  and
control their dependence with respect to $r$.
\bigskip

Let us define $\varphi(\gamma,\al)$ such that
$\varphi(\gamma,\al)= 2 \al $ if $\gamma=1$, $\varphi(\gamma,\al)
> 1 + \frac{2 \al-1}{\gamma} $ if $\frac{1-2\al}{1-\al} \le
\gamma<1$  and $\varphi(\gamma,\al) = \al$ if $0 \le \gamma <
\frac{1-2\al}{1-\al}$. Note that $\varphi(\gamma,\al) \in
[\al,2\al]$.

\bigskip

\begin{prop}\label{desh3}
Assume that $\si$ satisfies {\bf (H1)} and {\bf (H3)}. If $f \in W_0^{\al,\infty}(-r,T;\R^d)$ then
$$ \Vert G^{(\si)}_r (f) \Vert_{\al,\la}  \le \Lambda_\al(g) d^{(5)}  \left( 1+  \frac{\Vert f \Vert_{\al,\la(r)}}{\la^{1-\varphi(\gamma,\al)}} \right)$$
for all $\la \ge 1$ where $d^{(5)}$ is a positive constant depending only on $\al,\beta, T, d, m$
and $B_{0\al} = \Vert b_0 \Vert_{L^{1/\al}}$.
\end{prop}

\begin{dem}
From (\ref{NR418}) we can write
\begin{eqnarray}\label{auxi1}
& &\vert G^{(\si)}_r (f) (t) \vert   \cr & & \quad \le
\Lambda_\al(g) \left( \int_0^t \frac{\vert \si(s,f(s-r))
\vert}{s^\al} ds  + \al \int_0^t \int_0^s \frac{ \vert
\si(s,f(s-r)) - \si(u,f(u-r)) \vert}{ (s-u)^{\al + 1} } du ds
\right) \cr & & \quad \le \Lambda_\al(g) \left( K_0 \int_0^t
\frac{1+\vert f(s-r)\vert^\gamma}{s^\al} ds  + \al M_0 \int_0^t
\int_0^s \frac{ \vert f(s-r) - f(u-r) \vert}{ (s-u)^{\al + 1} } du
ds + C_\alpha \right)\cr & & \quad \le \Lambda_\al(g) C_\al \left(
1+\int_{-r}^{t-r} \frac{\vert f(s)\vert}{(s+r)^\al} ds  +
\int_{-r}^{t-r} \int_{-r}^{s} \frac{ \vert f(s) - f(u) \vert}{
(s-u)^{\al + 1} } du ds \right),
\end{eqnarray}
where we have used that $\vert f(s)\vert^\gamma \le 1 + \vert
f(s)\vert.$ 
On the other hand, using (\ref{NR417}) we have \goodbreak
\begin{eqnarray}\label{auxi2}
& & \int_0^t \frac{ \vert G^{(\si)}_r (f) (t) - G^{(\si)}_r (f)
(s) \vert}{(t-s)^{\al+1}} ds\cr & & \le \Lambda_\al(g) \left(
C_\al \int_0^t \frac{ \vert \si(s, f(s-r))\vert}{(t-s)^{2\al}} ds
+ \int_0^t \int_0^s \frac{ \vert \si(s,f(s-r)) -\si(u, f(u-r))
\vert}{ (s-u)^{\al + 1} } (t-u)^{-\al} du ds \right)\cr & &  \le
C_\al \Lambda_\al(g) \left( \int_{0}^{t} \frac{1+\vert
f(s-r)\vert^\gamma}{(t-s)^{ 2\al}} ds + \int_{0}^{t}
\frac{1}{(t-s)^\al} \int_{0}^{s} \frac{ \vert f(s-r) - f(u-r)
\vert}{ (s-u)^{\al + 1} } du ds  + C_\al \right) \cr & & \le C_\al
\Lambda_\al(g) \left( 1+ \int_{-r}^{t-r} \frac{\vert
f(s)\vert^\gamma}{(t-r-s)^{ 2\al}} ds  + \int_{-r}^{t-r}
\frac{1}{(t-r-s)^\al} \int_{-r}^{s} \frac{ \vert f(s) - f(u)
\vert}{ (s-u)^{\al + 1} } du ds  \right).
\end{eqnarray}

From H\"older inequality we can get that
\begin{eqnarray}\label{surtvarphi}
\int_{-r}^{t-r} \frac{\vert f(s)\vert^\gamma}{(t-r-s)^{ 2\al}} ds
& \le &
 C_\al \left( \int_{-r}^{t-r} \frac{\vert f(s)\vert}{(t-r-s)^{ \varphi(\gamma,\al)}} ds \right)^\gamma t^{1-\gamma-2\alpha +
 \varphi(\gamma,\al) \gamma}\\
& \le &
  C_\alpha \left( 1+ \int_{-r}^{t-r} \frac{\vert f(s)\vert}{(t-r-s)^{ \varphi(\gamma,\al)}} ds
 \right).\nonumber
\end{eqnarray}

Since $\al <\varphi(\gamma,\al)\le 2 \al$ and  $1-\gamma-2\alpha +
 \varphi(\gamma,\al) \gamma \ge 0$,
and
 and putting together (\ref{auxi1}) and
(\ref{auxi2}), we get that
$$
\Vert G^{(\si)}_r (f) \Vert_{\al,\la}  \le  C_\al \Lambda_\al(g)
\left( 1+  \left( \int_{-r}^{t-r}
\frac{e^{-\lambda(t-s)}}{(s+r)^\al} ds + \int_{-r}^{t-r}
\frac{e^{-\lambda(t-s)}}{(t-r-s)^{\varphi(\gamma,\al)}} ds \right)
\Vert f \Vert_{\al,\la(r)} \right).$$
We finish the proof using the fact that
$$ \int_0^t \frac{e^{-\lambda(t-s)}}{(t-s)^\al} \le \lambda^{\al-1}
\Gamma(1-\al) \quad {\rm and} \quad \int_0^t
\frac{e^{-\lambda(t-s)}}{s^\al} \le C_\al \lambda^{\al-1}.$$

 \hfill $\Box$
\end{dem}

\begin{obs}
Proposition \ref{desh3} is also true when $r=0$.
\end{obs}

\section{ Deterministic integral equations}

In this section we will study the deterministic delay  equations.
Following the method presented in \cite{NR}, we will prove a result of
existence and uniqueness of solution. We will also obtain a bound for
the $\Vert . \Vert_{\al,\la(r)}$ norm of the solution. In order to
obtain a bound whose dependence in $r$ could be controlled, we will
introduce a new method to compute this estimate.

\bigskip

Set $0 < \al < 1/2$ , $g \in W_T^{1-\al,\infty}(0,T;\R^m)$ and $\eta \in
W_0^{\al,\infty}(-r,0;\R^d) \cap C^{1-\al} (-r,0; \R^d).$ Consider the deterministic stochastic
differential equation on $\R^d$
\begin{eqnarray}\label{equag}
x (t) &= &\eta (0) + \int_0^t b(s, x) ds + \int_0^t \si(s, x(s-r)) dg_s, \quad t \in (0,T],\cr
x (t) &= &\eta (t), \qquad t \in [-r,0].
\end{eqnarray}
Using the notations introduced in the previous sections we can give another expression for equation
(\ref{equag}):
\begin{eqnarray*}
x (t) &= &\eta (0) + F^{(b)}(x)(t) + G_r^{(\si)} (x) (t), \quad t \in (0,T],\cr
x (t) &= &\eta (t), \qquad t \in [-r,0].
\end{eqnarray*}

\medskip

The result of existence and uniqueness reads as follows.

 \bigskip

\begin{teo}\label{teoexi}
Assume that $b$ and $\si$ satisfy hypothesis {\bf (H1)} and  {\bf (H2)} with
$\rho=1/\al, 0 <\beta,\delta\le 1$ and $$0 < \al < \al_0:=\min \{\frac12, \beta, \frac{\delta}{1+\delta} \}.$$
Then equation (\ref{equag}) has an unique solution $x \in W_0^{\al,\infty} (-r,T; \R^d) \cap
C^{1-\al}(-r,T;\R^d)$.
\end{teo}

\begin{dem}

{\it Step 1: $x \in C^{1-\al}(-r,T;\R^d).$}

If $x \in W_0^{\al,\infty} (-r,T; \R^d) $ is a solution then
$G_r^{(\si)} (x) \in C^{1-\al}(0,T;\R^d)$ (see Proposition
\ref{desint}) and $F^{(b)} (x) \in C^{1-\al}(0,T;\R^d)$ (see
Proposition \ref{desleb}). Furthermore
\begin{eqnarray*}
\Vert x \Vert_{1-\al(r)} & = & \Vert x \Vert_{\infty(r)} + \sup_{-r
\le s < t \le T}  \frac{\vert x(t)-x(s) \vert}{ (t-s)^{1-\al}} \cr
&\le & \Vert \eta \Vert_{\infty(-r,0)} +  \Vert F^{(b)} (x) \Vert_{\infty} + \Vert G_r^{(\si)} (x)  \Vert_{\infty}
 + \sup_{0
\le s < t \le T}  \frac{\vert x(t)-x(s) \vert}{ (t-s)^{1-\al}} \cr
& & \qquad + \sup_{-r \le s < t \le 0}  \frac{\vert
\eta(t)-\eta(s) \vert}{ (t-s)^{1-\al}} + \sup_{-r \le s \le 0 \le
t \le T} \frac{\vert x(t)-\eta(s) \vert}{ (t-s)^{1-\al}}
 \cr
&\le & \Vert \eta \Vert_{1-\al(-r,0)} +  \Vert F^{(b)} (x) \Vert_{1-\al} + \Vert G_r^{(\si)} (x)  \Vert_{1-\al}
 \cr
&& \qquad \quad + \sup_{-r \le s \le 0 \le t \le T} \left(
\frac{\vert x(t)-x(0) \vert}{ t^{1-\al}}+ \frac{\vert x(0)-\eta(s)
\vert}{ (-s)^{1-\al}} \right)
 \cr
&\le & 2 (\Vert \eta \Vert_{1-\al(-r,0)} +  \Vert F^{(b)} (x) \Vert_{1-\al} + \Vert G_r^{(\si)} (x)  \Vert_{1-\al})
< \infty.
\end{eqnarray*}

\goodbreak

{\it Step 2: Uniqueness}

 Consider $x$ and $x'$ two solutions such that $\Vert x \Vert_{1-\al(r)}\le N$ and
$\Vert x' \Vert_{1-\al(r)} \le N.$

Note that
$$
\sup_{t \in [-r,T]} e^{-\la t} \vert x(t)- x'(t) \vert_ \le
\sup_{t \in [0,T]} e^{-\la t} \vert F^{(b)} (x)(t) - F^{(b)}
(x')(t) \vert + \sup_{t \in [0,T]} e^{-\la t} \vert G_r^{(\si)}
(x)(t)- G_r^{(\si)} (x')(t) \vert
$$
and
\begin{eqnarray*}
&&\sup_{t \in [-r,T]} e^{-\la t} \int_{-r}^t \frac{\vert x(t)-x'(t)-(x(s)-x'(s)) \vert}{(t-s)^{\al+1}} ds
\cr
&& \quad  =
\sup_{t \in [0,T]} e^{-\la t} \int_{-r}^t \frac{\vert x(t)-x'(t)-(x(s)-x'(s)) \vert}{(t-s)^{\al+1}} ds \cr
&& \quad \le \sup_{t \in [0,T]} e^{-\la t} \int_{0}^t \frac{\vert x(t)-x'(t)-(x(s)-x'(s)) \vert}{(t-s)^{\al+1}} ds
+ \sup_{t \in [0,T]} e^{-\la t} \int_{-r}^0 \frac{\vert x(t)-x'(t) \vert}{(t-s)^{\al+1}} ds.
\end{eqnarray*}
So
\begin{equation}\label{xx}
\Vert x- x' \Vert_{\al,\la(r)} \le  \Vert F^{(b)} (x) - F^{(b)} (x')  \Vert_{\al,\la} + \Vert  G_r^{(\si)} (x)-  G_r^{(\si)} (x') \Vert_{\al,\la} + U,
\end{equation}
where
\begin{equation*}
U:=\sup_{t \in [0,T]} e^{-\la t} \int_{-r}^0 \frac{\vert x(t)-x'(t) \vert}{(t-s)^{\al+1}} ds.
\end{equation*}

Let us study $U$. Clearly,
\begin{equation}\label{desu}
U \le \sup_{t \in [0,T]} e^{-\la t} \vert x(t)-x'(t) \vert  \frac{1}{\al t^{\al}} \le U_1+U_2,
\end{equation}
with
\begin{equation*}
U_1:=  \sup_{t \in [0,T]}  \frac{e^{-\la t}}{\al t^{\al}}  \vert F^{(b)} (x)(t) - F^{(b)} (x')(t) \vert ,\quad
U_2:=  \sup_{t \in [0,T]}  \frac{e^{-\la t}}{\al t^{\al}}  \vert G_r^{(\si)} (x)(t) - G_r^{(\si)} (x')(t) \vert.
\end{equation*}

Moreover
\begin{eqnarray}\label{desu1}
U_1 & \le &  \frac{1}{\al} \sup_{t \in [0,T]}  \frac{e^{-\la
t}}{t^{\al}}  \vert \int_0^t (b(s,x)-b(s,x'))ds \vert \le
\frac{L_N}{\al} \sup_{t \in [0,T]}  \frac{e^{-\la t}}{t^{\al}}
\int_0^t  \sup_{-r \le u \le s} \vert x(u)- x'(u) \vert ds\cr &
\le &  \frac{L_N}{\al} \left( \sup_{u \in [-r,T]}  e^{-\la u}
\vert x(u)- x'(u) \vert \right)
 \sup_{t \in [0,T]} \int_0^t \frac{e^{-\la (t-s)}}{(t-s)^{\al}}   ds \cr
& \le &
 \frac{L_N}{\al} \la^{\al-1} \Gamma(1-\al) \Vert x- x'
 \Vert_{\al,\la(r)}.
 \end{eqnarray}

On the other hand, using (\ref{NR418}), we can write
\begin{eqnarray*}
U_2 & \le &  \frac{1}{\al} \sup_{t \in [0,T]}  \frac{e^{-\la t}}{t^{\al}}  \vert \int_0^t (\sigma(s,x(s-r))-\sigma(s,x'(s-r)))dg_s \vert \cr
&\le&  \frac{1}{\al} \sup_{t \in [0,T]}  \frac{e^{-\la t}}{t^{\al}} \Lambda_\al (g) \left(  \int_0^t
 \frac{\vert \sigma(s,x(s-r))-\sigma(s,x'(s-r)) \vert}{s^\al} ds  \right. \cr
  & & \qquad \left. + \al
 \int_0^t \int_0^s
 \frac{\vert \sigma(s,x(s-r))-\sigma(s,x'(s-r)) - \sigma(u,x(u-r))+ \sigma(u,x'(u-r))\vert}{(s-u)^{\al+1}} du ds \right)
 \cr
&\le&   C_\al \Lambda_\al (g) \sup_{t \in [0,T]}  \left(
\int_0^t \frac{e^{-\la (t-s)}}{s^{2\al}} ds\right)
  \sup_{s \in [0,T]}   \
  \Big[ e^{-\la s} (\vert \sigma(s,x(s-r))-\sigma(s,x'(s-r)) \vert \Big. \cr
  & & \qquad \Big. +
 e^{-\la s} \int_0^s
 \frac{\vert \sigma(s,x(s-r))-\sigma(s,x'(s-r)) - \sigma(u,u(y-r))+ \sigma(u,x'(u-r))\vert}{(s-u)^{\al+1}} du \Big].
 \end{eqnarray*}
 Now, using Lemma \ref{calcDavid} and that
 $$\int_0^t \frac{e^{-\la (t-s)}}{s^{2\al}} ds
 \le C_\al \la^{2 \al -1}$$ we obtain
 \begin{equation}\label{desu2}
 U_2 \le C_\al \la^{2 \al -1} \Lambda_\al (g) ( 1 + \Delta_r(x) + \Delta_r (x') )
 \Vert x - x' \Vert_{\al,\la(r)}.
 \end{equation}

Putting together (\ref{xx}), (\ref{desu}), (\ref{desu1}) and
(\ref{desu2}) and using Propositions \ref{desleb} and \ref{desint}
we get that for all $\lambda \ge 1$
$$
\Vert x - x' \Vert_{\al,\la(r)} \le C_\al C_N^{(1)} \left(
\frac{1}{\la^{1-\al}} + \frac{1}{\la^{1-2\al}}  \Lambda_\al (g) (
1 + \Delta_r(x) + \Delta_r (x') ) \right)
 \Vert x - x' \Vert_{\al,\la(r)}.
 $$

Finally, since
\begin{eqnarray*}
& &\Delta_r(x) = \sup_{s \in [-r,T]} \int_{-r}^s \frac{ \vert
x(s)-x(u) \vert^\delta}{(s-u)^{\al+1}} ds \le N \sup_{s \in
[-r,T]} \int_{-r}^s \frac{ (s -u)^{(1-\al)\delta}}{(s-u)^{\al+1}}
ds \cr & & \qquad = \frac{
N(T+r)^{\delta-\al(1+\delta)}}{\delta-\al(1+\delta)}:=C_N^{(2)}
\end{eqnarray*}
choosing $\la$ large enough such that
$$
C_\al C_N^{(1)} \left( \frac{1}{\la^{1-\al}} +
\frac{1}{\la^{1-2\al}}  \Lambda_\al (g) ( 1 + 2C_N^{(2)} ) \right)
\le \frac12$$ we get that
$$\frac12 \Vert x - x' \Vert_{\al,\la(r)} \le 0$$
and obviously $x=x'$.

\medskip

{\it Step 3: Existence}

Let us consider the operator $\cL: W_0^{\al,\infty} (-r,T;\R^d)
\to C^{1-\al}(-r,T;\R^d)$ such that
\begin{eqnarray}\label{equaope}
\cL(y)(t) &= &\eta (0) + \int_0^t b (s, y) ds + \int_0^t \si (s,
y(s-r)) dg_s, \quad t \in (0,T],\cr \cL(y)(t) &= &\eta (t), \qquad
t \in [-r,0].
\end{eqnarray}
Let us use the notatiton $y^*=\cL(y)$

In order to prove the existence of $y$ such that $y=\cL(y)$ we will
use a fixed point argument based in Lemma \ref{puntfixe}. We will
check the three conditions of that lemma.

\smallskip

{\sl Condition 1.} Note that for $t \in [-r,0]$
$$
\vert y^*(t)\vert + \int_{-r}^t \frac{\vert y^*(t)-y^*(s) \vert}{
(t-s)^{\alpha+1}}ds=  \vert \eta(t)\vert + \int_{-r}^t \frac{\vert
\eta(t)-\eta(s) \vert}{ (t-s)^{\alpha+1}}ds
$$
and for $t \in (0,T]$
$$
\vert y^*(t)\vert + \int_{-r}^t \frac{\vert y^*(t)-y^*(s) \vert}{
(t-s)^{\alpha+1}}ds= \vert y^*(t)\vert + \int_{-r}^0 \frac{\vert
y^*(t)-\eta(s) \vert}{ (t-s)^{\alpha+1}}ds+ \int_{0}^t \frac{\vert
y^*(t)-y^*(s) \vert}{ (t-s)^{\alpha+1}}ds.
$$
Hence
\begin{equation}\label{exis1}
\Vert y^* \Vert_{\al,\la(r)} \le \Vert \eta \Vert_{\al,\la(-r,0)}+
\Vert F^{(b)}(y) \Vert_{\al,\la}+\Vert G_r^{(\si)}(y)
\Vert_{\al,\la} + E
\end{equation}
where
$$
E:= \sup_{t \in [0,T]} e^{-\la t} \int_{-r}^0 \frac{\vert
y^*(t)-\eta(s) \vert}{ (t-s)^{\alpha+1}}ds.$$

Let us study $E$. Clearly $E \le E_1 + E_2$ with
\begin{eqnarray}\label{exis2}
& &E_1:= \sup_{t \in [0,T]} e^{-\la t} \int_{-r}^0 \frac{\vert
y^*(t)-\eta(0) \vert}{ (t-s)^{\alpha+1}}ds \le C_\alpha \sup_{t
\in [0,T]} \frac{e^{-\la t}}{t^\al} \vert y^*(t)-\eta(0) \vert,\cr
& &E_2:= \sup_{t \in [0,T]} e^{-\la t} \int_{-r}^0 \frac{\vert
\eta(0)-\eta(s) \vert}{ (-s)^{\alpha+1}}ds \le \Vert \eta
\Vert_{\al,\la(-r,0)}\label{exis3}.
\end{eqnarray}
In order to study $E_1$, we shall repeat similar computations to
those used to estimate $U$ when we proved the uniqueness. We will
give only a sketch of these computations. Observe first that $E_1
\le E_{1,1}+ E_{1,2}$ where
$$
E_{1,1}:= C_\alpha \sup_{t \in [0,T]} \frac{e^{-\la t}}{t^\al}
\int_0^t \vert b(s,y) \vert ds, \qquad E_{1,2}:= C_\alpha \sup_{t
\in [0,T]} \frac{e^{-\la t}}{t^\al} \vert \int_0^t \si(s,y(s-r))
dg_s \vert.
$$
First
\begin{eqnarray}\label{exis4}
E_{1,1} & \le & C_\al L_0 \left( \sup_{s \in [-r,T]} e^{-\la s}
\vert y(s)\vert \right)  \sup_{t \in [0,T]}
 \int_0^t \frac{e^{-\la (t-s)}}{(t-s)^\al} ds +  C_\al \sup_{t \in [0,T]} \frac{e^{-\la t}}{t^\al}
  \int_{0}^t b_0(s) ds \cr
  & \le & C_\al L_0 \la^{\al-1} \Gamma(\al-1) \Vert y \Vert_{\al,\la(r)} + C_\al  B_{0,\al}.
\end{eqnarray}
On the other hand, using (\ref{NR418})
\begin{eqnarray}\label{exis5}
E_{1,2} & \le & C_\al \sup_{t \in [0,T]}  \frac{e^{-\la t}}{t^\al}
\Lambda_\al(g) \left(
 \int_0^t \frac{M_0(s^\beta + \vert y(s-r)\vert)}{s^\al} ds \right. \cr
 & & \qquad \qquad  \left. +
 \al \int_0^t \int_0^s M_0 \frac{( \vert s-u\vert^\beta  +
 \vert y(s-r)-y(u-r)\vert)}{(s-u)^{\al+1}} duds \right) \cr
  & \le & C_\al \Lambda_\al(g)  \left(e^{-\la r} \Vert y \Vert_{\al,\la(r)}
  \sup_{t \in [0,T]} \left( \int_0^t \frac{e^{-\la (t-s)}}{s^{2\al}} ds +
 \int_0^t \frac{e^{-\la (t-s)}}{(t-s)^\al} ds\right) \right. \cr
 & & \qquad \qquad\left. +
  \sup_{t \in [0,T]} \frac{e^{-\la t}}{t^\al}  \int_0^t s^{-\al}  ds +
   \sup_{t \in [0,T]} \frac{e^{-\la t}}{t^\al}  \int_0^t \int_0^s  (s-u)^{\beta-\al-1} du  ds \right)\cr
    & \le &  \Lambda_\al(g) C_\al ( \la^{\al-1} \Gamma(\al-1) + \la^{2\al-1} ) \Vert y \Vert_{\al,\la(r)}
    e^{-\la r}+  \Lambda_\al(g) C_\al T^{1+\beta-2\al}.
\end{eqnarray}
Now putting together (\ref{exis1})-(\ref{exis5}) and using
Propositions \ref{desleb} and \ref{desint} in order to bound the
terms $\Vert F^{(b)}(y) \Vert_{\al,\la}$ and $\Vert G_r^{(\si)}(y)
\Vert_{\al,\la}$ respectively, we obtain that
\begin{equation}\label{exis6}
\Vert y^* \Vert_{\al,\la(r)} \le M_1 (\la) + M_2 (\la) (1+  \Vert
y \Vert_{\al,\la(r)})  +  M_3 (\la) \Vert y \Vert_{\al,\la(r)}
\end{equation}
where
\begin{eqnarray}\label{exist}
&&M_1(\la) := 2\Vert \eta \Vert_{\al,\la(-r,0)}+  \Lambda_\al(g)
C_\al T^{1+\beta-2\al} + C_\al  B_{0,\al},\qquad M_2(\la):=
(\frac{d^{(2)}}{\la^{1-2\al}} + \frac{\Lambda_\al(g)
d^{(4)}}{\la^{1-2\al}} ) \cr && \qquad  M_3(\la) := C_\al \left(
\frac{L_0}{\la^{1-\al}} \Gamma(\al-1) + \Lambda_\al(g)  (
\frac{\Gamma(\al-1)}{\la^{1-\al}}  + \frac{1}{\la^{1-2\al}} )
e^{-\la r}\right).
\end{eqnarray}
Choose $\la=\la_0$ large enough in order to be $M_2(\la_0) +
M_3(\la_0) \le \frac12$. If $\Vert y \Vert_{\al,\la_0(r)} \le 2 ( 1+
M_1(\la_0))$ then $ \Vert y^* \Vert_{\al,\la_0(r)} \le 2 ( 1+
M_1(\la_0))$ and so $\cL (B_0) \subset B_0$ where
$$
B_0 := \{ y \in W_0^{\al,\infty} (-r,T;\R^d), \Vert y
\Vert_{\al,\la_0(r)} \le 2 ( 1+ M_1(\la_0)) \}.$$ So, {\sl Condition
1} is satisfied with the metric $\rho_0$ associated to the norm $\Vert
\cdot \Vert_{\al,\la_0(r)}.$

\medskip

{\sl Condition 2.} Notice first that if $y \in B_0$ then $\Vert y \Vert_{\al,\infty(r)} \le 2 e^{\lambda_0 T}
(1 + M_1(\lambda_0)):=N_0$. Then, repeating the same computations we have done in
step 2 we get that for all $y,y' \in B_0$ and for all $\la \ge 1$
\begin{eqnarray}\label{cond2}
\Vert \cL(y) - \cL(y') \Vert_{\al,\la(r)} & \le & C_\al C_{N_0}
\left( \frac{1}{\la^{1-\al}} + \frac{1}{ \la^{1-2 \al } }
\Lambda_\al (g) ( 1 + \Delta_r(y) + \Delta_r (y') ) \right)
 \Vert y - y' \Vert_{\al,\la(r)}\cr
 & \le & \frac{1}{\la^{1-2 \al}}
C_1( 1 + \Delta_r(y) + \Delta_r (y') )
 \Vert y - y' \Vert_{\al,\la(r)},
\end{eqnarray}
with $ C_1:=C_\al C_{N_0} \Lambda_\al (g)$ and where we recall
that
\begin{equation*}
\Delta_r(y) = \sup_{s \in [-r,T]} \int_{-r}^s \frac{ \vert
y(s)-y(u) \vert^\delta}{(s-u)^{\al+1}} ds.
\end{equation*}
Note that $\Delta_r: W_0^{\al,\infty} (-r,T;\R^d) \to [0,+\infty]$ is a lower semicontinuous function.

Given $y \in \cL (B_0)$ consider $\hat y \in B_0$ such that $y=\cL
(\hat y) \in C^{1-\al} (-r,T; \R^d)$. Repeating the computations of
step 1 and using Propositions \ref{desint} and \ref{desleb} we get
\begin{eqnarray*}
\Vert y \Vert_{1-\al(r)} &\le & 2  (\Vert \eta \Vert_{1-\al(-r,0)}
+ (d^{(1)}+ \Lambda_\al (g) d^{(3)} )  ( 1+ \Vert \hat y
\Vert_{\al,\infty(r)})) \cr &\le &  2 (d^{(1)}+ \Lambda_\al (g)
d^{(3)} )  ( 1+
 2 e^{\lambda_0 T}
(1 + M_1(\lambda_0)))
+ 2  \Vert \eta
\Vert_{1-\al(-r,0)} := C_2.
\end{eqnarray*}
Hence
\begin{eqnarray*}
& &\Delta_r(y)=\sup_{-r \le u \le T} \int_{-r}^u \frac{ \vert y (u)-
y (s) \vert^\delta}{(u-s)^{\al+1}} ds  \le \sup_{-r \le u \le T}
\int_{-r}^u \frac{ \vert (u-s)^{1-\al} \vert^\delta \Vert y
\Vert_{1-\al(r)}^\delta }{(u-s)^{\al+1}} ds\\
& & \qquad \le \frac{(T+r)^{\delta(1-\al)-\al}}{\delta(1-\al)-\al}:= C_3.
\end{eqnarray*}
So {\sl Condition 2} is fullfilled for the metric associated with the
norm $\Vert \cdot \Vert_{\al,1(r)}$ and $\varphi(y)=C_1 (\frac12 +
\Delta_r(y)), C_0=C_1 (\frac12 + C_3)$ and $K_0=2C_0$.

\medskip

{\sl Condition 3.}  From the computations in the proof of the
above condition we get that for all $y, y' \in  \cL (B_0)$
$$
\Vert \cL (y) - \cL (y') \Vert_{\al,\la(r)} \le C_1 \frac{1 + 2
C_3}{\la^{1-2\al}} \Vert y-y' \Vert_{\al,\la(r)}.
$$
So, it suffices to choose $\la=\la_2$ such that
$$
C_1 \frac{1 + 2 C_3}{\la_2^{1-2\al}} \le \frac12. $$

\hfill $\Box$
\end{dem}

\bigskip

We finish this section providing an upper bound for the norm of the
solution. We obtain the bound as a consequence of the previous
computations and without making use of Gronwall's lemmas.

\smallskip

Let us recall the definition of $\varphi(\gamma,\al)$. Set
$\varphi(\gamma,\al)= 2 \al $ if $\gamma=1$, $\varphi(\gamma,\al)
> 1 + \frac{2 \al-1}{\gamma} $ if $\frac{1-2\al}{1-\al} \le
\gamma<1$  and $\varphi(\gamma,\al) = \al$ if $0 \le \gamma <
\frac{1-2\al}{1-\al}$.

\begin{lema}\label{supnorminfprim} Assume {\bf (H1)}, {\bf (H2)} and {\bf (H3)}.
Then the unique solution of equation  (\ref{equag})  satisfies
$$
\Vert x \Vert_{\al,\infty(r)} \le d_\al^{(6)} (\Vert \eta
\Vert_{\al,\infty(-r,0)}+  \Lambda_\al(g) + 1) e^{d_\al^{(7)} +
d_\al^{(8)} \Lambda_\al(g)^{\frac{1}{1-\varphi(\al,\gamma)}}}.
$$
\end{lema}

\begin{dem}
We have to repeat the same computations we have done in  the proof of
Theorem \ref{teoexi} ({\it Step 3, Condition 1}). Note only that now
we have hypothesis  {\bf (H3)}, so we will make use of Proposition
\ref{desh3} and we will repeat the computations of the term $E_{1,2}$
using  {\bf (H3)} and the ideas of the proof of Proposition
\ref{desh3}. Then, we obtain that
\begin{equation*}
\Vert x \Vert_{\al,\la(r)} \le M_{1} (\la) + M_{2} (\la)  \Vert x
\Vert_{\al,\la(r)}
\end{equation*}
where
\begin{eqnarray*}
M_{1} (\la) := 2\Vert \eta \Vert_{\al,\la(-r,0)}+ \Lambda_\al(g)
C_\al  + C_\al B_{0,\al}+ \Lambda_\al(g) d^{(5)} +
d^{(2)} ,
 \cr
  M_2(\la) := C_\al \left( \frac{L_0}{\la^{1-\al}} + \Lambda_\al(g)   \frac{
d^{(5)}+1}{\la^{1-\varphi(\gamma,\al)}}  +
\frac{d^{(2)}+1}{\la^{1-\al}} \right).
\end{eqnarray*}
Choose $\la=\la_0$ large enough such that $M_2(\la_0) \le \frac12$.
Then
$$
\Vert x \Vert_{\al,\la_0(r)} \le 2 \left( 2\Vert \eta
\Vert_{\al,\la_0(-r,0)}+  \Lambda_\al(g) C_\al  +
C_\al  B_{0,\al} + \Lambda_\al(g) d^{(5)}  + d^{(2)} \right).
$$
Note that
\begin{eqnarray*}
\la_0 & \le  &  2 C_\al \left(  L_0 + \Lambda_\al(g)  ( d^{(5)} + 1) + d^{(2)} + 1 \right)^{\frac{1}{1-\varphi(\gamma,\al)}}\\
& \le & d_\al \left( 2 C_\al (L_0+  d^{(2)} + 1) \right)^{\frac{1}{1-\varphi(\gamma,\al)}}  +
\Lambda_\al(g)^{\frac{1}{1-\varphi(\gamma,\al)}} d_\al \left(2 C_\al ( 1+ d^{(5)} ) \right)^{\frac{1}{1-\varphi(\gamma,\al)}}.
\end{eqnarray*}
Hence
$$
\Vert x \Vert_{\al,\infty(r)} \le K_\al
e^{d_\al(1+\Lambda_\al(g)^{\frac{1}{1-\varphi(\gamma,\al)}})},
$$
with $$ K_\al = 2 \left( 2\Vert \eta \Vert_{\al,\infty(-r,0)}+
\Lambda_\al(g) C_\al  + C_\al B_{0,\al}+
\Lambda_\al(g) d^{(5)}  + d^{(2)} \right)$$ and the proof finishes
easily.

 \hfill $\Box$
\end{dem}

\begin{obs}
Note that {\bf (H1)} implies  {\bf (H3)} with $\gamma=1$.
\end{obs}

\section {Convergence when the delay goes to zero}

Our aim here is to study what happens when the delay $r$ tends to
zero. We will assume the hypothesis {\bf (H1)} and {\bf (H2')}
throughout this section. Observe that all the results given under assumption
 {\bf (H2)} in the previous sections also hold under assumption
 {\bf (H2')}.

 Set $x^r$ the solution of the integral delay
equation on $\R^d$ \goodbreak
\begin{eqnarray}\label{equar}
x^r(t) &= &\eta (0) + \int_0^t b (s, x^r(s)) ds + \int_0^t \si (s, x^r(s-r)) dg_s, \quad t \in (0,T],\cr
x^{r} (t) &= &\eta (t), \qquad t \in [-r,0],
\end{eqnarray}
and $x$ the solution of the integral equation
on $\R^d$
\begin{equation}\label{equasenser}
x(t) = \eta (0) + \int_0^t b (s, x(s)) ds + \int_0^t \si (s, x(s)) dg_s, \quad t \in (0,T].
\end{equation}
From the previous sections and the paper of Nualart and Rascanu
\cite{NR}, we know that these solutions exist, they are unique and
$x^r \in W_0^{\al,\infty} (-r,T;\R^d)$ and $x \in W_0^{\al,\infty}
(0,T;\R^d)$

\medskip

Let us start by proving two technical lemmas that we will use in the
sequel.

\bigskip

\begin{lema}\label{supnorminf} Assume {\bf (H1)} and {\bf (H2')}.
Suppose that there exists $r_0 >0$ such that $$\eta \in
W_0^{\al,\infty}(-r_0,0;\R^d)  \cap  C^{1-\al} (-r_0,0; \R^d).$$
Then
$$
\Vert x^r \Vert_{\al,\infty(r)} \le C_\al^{(0)} (\Vert \eta
\Vert_{\al,\infty(-r,0)}+  \Lambda_\al(g) + 1) e^{C_\al^{(1)} +
C_\al^{(2)} \Lambda_\al(g)^{\frac{1}{1-2\al}}}
$$
and
$$
\sup_{0 \le r \le r_0} \Vert x^r \Vert_{\al,\infty(r)} \le
C_\al^{(3)}.
$$
Moreover, if {\bf (H3)} holds, then
$$
\Vert x^r \Vert_{\al,\infty(r)} \le C_\al^{(0)} (\Vert \eta
\Vert_{\al,\infty(-r,0)}+  \Lambda_\al(g) + 1) e^{C_\al^{(1)} +
C_\al^{(2)} \Lambda_\al(g)^{\frac{1}{1-\varphi(\gamma,\alpha)}}}
$$
\end{lema}

\begin{dem}
Follow the  ideas  of Proposition
\ref{supnorminfprim} using also that
$\sup_{0 \le r \le r_0} \Vert \eta \Vert_{\al,\infty(-r,0)} < \infty.$

 \hfill $\Box$
\end{dem}

\medskip

\begin{lema}\label{supdeltar}
Assume {\bf (H1)} and {\bf (H2')}. Suppose that there exists $r_0 >0$ such that $$\eta \in
W_0^{\al,\infty}(-r_0,0;\R^d) \cap C^{1-\al} (-r_0,0; \R^d)$$  and  that $\al  <
\delta/(1+\delta)$ then
$$
\sup_{0 \le r \le r_0} \Delta_r(x^r) \le C_\al^{(4)}.
$$
\end{lema}

\begin{dem}
From the definition of $\Delta_r$ we have
\begin{eqnarray}\label{aa}
& &\Delta_r(x^r)=\sup_{-r \le u \le T} \int_{-r}^u \frac{ \vert
x^r(u)-x^r(s) \vert^\delta}{(u-s)^{\al+1}} ds  \le \sup_{-r \le u
\le T} \int_{-r}^u \frac{ \vert (u-s)^{1-\al} \vert^\delta \Vert
x^r \Vert_{1-\al(r)}^\delta }{(u-s)^{\al+1}} ds \cr & & \qquad
\le C_\al \sup_{0 \le r \le r_0} \Vert x^r \Vert_{1-\al(r)}^\delta
 \sup_{-r \le u \le T} \int_{-r}^u \frac{1}{(u-s)^{\al+1-\delta +\delta \al}}
 ds.
\end{eqnarray} From Propositions \ref{desint} and \ref{desleb}
(see also step 1 in the proof of Theorem \ref{teoexi}) we get
\begin{eqnarray}\label{bb}
\Vert x^r \Vert_{1-\al(r)} &\le & 2 (\Vert \eta
\Vert_{1-\al(-r,0)} +  \Vert F^{(b)} (x^r) \Vert_{1-\al} + \Vert
G_r^{(\si)} (x^r) \Vert_{1-\al})\cr &\le & 2 (\Vert \eta
\Vert_{1-\al(-r,0)} + (d^{(1)}+ \Lambda_\al (g) d^{(3)} )  ( 1+
\Vert x^r \Vert_{\al,\infty(r)})) .
\end{eqnarray}
The result follows from (\ref{aa}), (\ref{bb}), Lemma
\ref{supnorminf} and the fact that $\delta- \al - \al \delta > 0$.
 \hfill $\Box$
\end{dem}

\bigskip

The main results of this section is the following theorem:

\begin{teo}\label{teoconver}
Assume that $b$ and $\si$ satisfy hypothesis {\bf (H1)} and  {\bf
(H2')} with $\rho=1/\al, 0 <\beta,\delta\le 1$ and $$0 < \al <
\al_0:=\min \{\frac12, \beta, \frac{\delta}{1+\delta} \}.$$ If there exits $r_0$ such that
$$\eta \in
W_0^{\al,\infty}(-r_0,0;\R^d)  \cap  C^{1-\al} (-r_0,0; \R^d)$$
 and $\al  < \delta/(1+\delta)$, then
$$
\lim_{r \to 0} \Vert x- x^r \Vert_{\al,\infty} = 0.
$$
\end{teo}

\begin{dem}
Actually, we will proof that there exists $\la_0$ such that
$$
\lim_{r \to 0} \Vert x- x^r \Vert_{\al,\la_0} = 0.
$$

Using Lemma \ref{supnorminf}, let us choose $N$ such that
$\Vert x \Vert_{\al,\infty} \le N$ and $\sup_{0 \le r \le r_0} \Vert
x^r \Vert_{\al,\infty(r)} \le N$. By Proposition \ref{desleb} we get
that
\begin{equation}\label{raobe}
\Vert F^{(b)}(x) -  F^{(b)}(x^r) \Vert_{\al,\la} \le \frac{d_N}{\la^{1-\al}} \Vert x  - x^r \Vert_{\al,\la}
\end{equation}
Let us define a function $y^r:[0,T] \to \R^d$ such that
$y^r(s):=x^r(s-r)$. Then $y^r \in  W_0^{\al,\infty} (0,T;\R^d)$ and it
is easy to check that $\Delta_0(y^r) \le \Delta_r (x^r)$ and
$G_r^{(\si)} (x^r) = G^{(\si)} (y^r)$. From Proposition \ref{desint}
whith $r=0$, we obtain that
\begin{equation}\label{raosigma}
\Vert G^{(\si)} (x) - G^{(\si)}_r (x^r) \Vert_{\al,\la}  \le
\frac{\Lambda_\al(g) d^{(2)}_N}{\la^{1-2\al}} (1 +
\Delta_0(x)+\Delta_r(x^r)) \Vert x-y^r \Vert_{\al,\la} .
\end{equation}

Thanks to (\ref{raobe}), (\ref{raosigma}) and the inequality
$$
\Vert x-y^r \Vert_{\al,\la} \le \Vert x-x^r \Vert_{\al,\la} + \Vert x^r-y^r \Vert_{\al,\la}
$$
we can write
\begin{eqnarray*}
\Vert x-x^r \Vert_{\al,\la} & \le & \left( \frac{d_N}{\la^{1-\al}}
+ \frac{\Lambda_\al(g) d^{(2)}_N}{\la^{1-2\al}} (1 +
\Delta_0(x)+\Delta_r(x^r)) \right) \Vert x  - x^r \Vert_{\al,\la}\cr
& & + \frac{\Lambda_\al(g) d^{(2)}_N}{\la^{1-2\al}} (1 +
\Delta_0(x)+\Delta_r(x^r)) \Vert x^r-y^r \Vert_{\al,\la}.
\end{eqnarray*}

Hence, if we choose $\la_0$ large enough such that
$$
\frac{d_N}{\la_0^{1-\al}} + \frac{\Lambda_\al(g)
d^{(2)}_N}{\la_0^{1-2\al}} (1 + \Delta(x)+ \sup_{-r_0 \le r <0}
\Delta_r(x^r)) \le \frac12 ,
$$
we obtain that for all $r \in (-r_0,0)$
$$
\Vert x-x^r \Vert_{\al,\la_0}  \le  \Vert x^r  - y^r
\Vert_{\al,\la_0}.$$
 So, to finish the proof it suffices to show
that
$$
\lim_{r \to 0} \Vert x^r- y^r \Vert_{\al,\la_0} = 0.
$$
Actually, we will check that
\begin{equation}\label{ultpas}
\lim_{r \to 0} \Vert x^r- y^r \Vert_{\al,\infty} = 0.
\end{equation}

Let us observe that
\begin{eqnarray*}
\vert x^r(t)-y^r(t) \vert &= & \vert x^r(t)-x^r(t-r) \vert \cr
&\le& \vert  F^{(b)}(x^r)(t) -  F^{(b)}(x^r)((t-r) \vee 0) \vert
+ \vert  G^{(\si)}_r (x^r)(t) -  G^{(\si)}_r (x^r)((t-r) \vee 0) \vert \cr
& & \qquad + \vert \eta(0) - \eta(0 \wedge (t-r)) \vert.
\end{eqnarray*}
From the fact that $F^{(b)}(x^r), G^{(\si)}_r (x^r)$ and $\eta$ are of $C^{1-\al}$ and using the estimates of the norms given in
Propositions \ref{desleb} and \ref{desint} we obtain easily that
\begin{equation}\label{raoholder}
\vert x^r(t)-y^r(t) \vert \le \left((d^{(1)} + \Lambda_\al(g)
d^{(3)} )(1+  \Vert x^r \Vert_{\al,\infty(r)}) + \Vert \eta
\Vert_{1-\al(-r,0)} \right) r^{1-\al}.
\end{equation}

On the other hand, we also have to deal with
\begin{equation*}
T_1:= \int_0^t \frac{\vert x^r(t) - x^r(t-r)-x^r(s)+x^r(s-r) \vert}{(t-s)^{\al + 1}} ds
\end{equation*}
Let us assume that $t>r$. (The case $t<r$ can be computed easily
following the ideas that we will use to study $T_{1,2}$) Clearly
$T_1:= T_{1,1} + T_{1,2} $ where
\begin{eqnarray*}
T_{1,1}&:=& \int_0^{t-r} \frac{\vert x^r(t) - x^r(t-r)\vert +\vert x^r(s)-x^r(s-r) \vert}{(t-s)^{\al + 1}} ds, \cr
T_{1,2}&:=& \int_{t-r}^t \frac{\vert x^r(t) - x^r(s)\vert + \vert x^r(t-r)-x^r(s-r) \vert}{(t-s)^{\al + 1}} ds,
\end{eqnarray*}
and
\begin{eqnarray*}
T_{1,1} & \le & 2 \left((d^{(1)} + \Lambda_\al(g) d^{(3)} )(1+
\Vert x^r \Vert_{\al,\infty(r)}) + \Vert \eta \Vert_{1-\al(-r,0)}
\right) r^{1-\al} \int_0^{t-r} \frac{1}{(t-s)^{\al + 1}} ds\cr & =
&  2 \left((d^{(1)} + \Lambda_\al(g) d^{(3)} )(1+  \Vert x^r
\Vert_{\al,\infty(r)}) + \Vert \eta \Vert_{1-\al(-r,0)} \right)
r^{1-\al}  \frac{1}{\alpha} \left( \frac{1}{r^\al} -
\frac{1}{t^\al} \right)\cr & \le  &  \frac{2}{\al} \left((d^{(1)}
+ \Lambda_\al(g) d^{(3)} )(1+  \Vert x^r \Vert_{\al,\infty(r)}) +
\Vert \eta \Vert_{1-\al(-r,0)} \right) r^{1-2\al},
\end{eqnarray*}
and
\begin{eqnarray*}
T_{1,2} & \le & 2 \left((d^{(1)} + \Lambda_\al(g) d^{(3)} )(1+
\Vert x^r \Vert_{\al,\infty(r)}) + \Vert \eta \Vert_{1-\al(-r,0)}
\right) \int_{t-r}^t \frac{(t-s)^{1-\al} }{(t-s)^{\al + 1}} ds\cr
& = &  \frac{2}{1-2\al} \left((d^{(1)} + \Lambda_\al(g) d^{(3)}
)(1+  \Vert x^r \Vert_{\al,\infty(r)}) + \Vert \eta
\Vert_{1-\al(-r,0)} \right) r^{1-2\al},
\end{eqnarray*}

So
\begin{equation}\label{dest1}
T_1 \le C_\al  \left((d^{(1)} + \Lambda_\al(g) d^{(3)} )(1+ \Vert
x^r \Vert_{\al,\infty(r)}) + \Vert \eta \Vert_{1-\al(-r,0)}
\right) r^{1-2\al}.
\end{equation}

Using Lemma \ref{supnorminf} and putting together
(\ref{raoholder}) and (\ref{dest1}), we get easily (\ref{ultpas})
and the proof is complete.

\hfill $\Box$
\end{dem}

\section{ Stochastic integral equations}

In this section we will apply the results of the previous two sections
in order to prove the main theorems of this paper.

\bigskip
The stochastic integral appearing throughout this paper $ \int_0^T
u(s) dW_s $ is a path-wise Riemann-Stieltjes integral and it is
well know that this integral exists if the process $u(s)$ has
H\"older continuous trajectories of order larger than $1-H$.

Set $\al \in (1-H, \frac12)$. For any $\delta \in (0,2)$, by
Fernique's theorem it holds that $$ E( \exp(\Lambda_\al (W)^\delta
)) < \infty.$$
Then if $u=\{u_t, t \in [0,T]\}$ is a stochastic
process whose trajectories belong to the space $W_T^{\al,1}
(0,T)$, the Riemann-Sieltjes integral $\int_0^T u(s) dW_s$ exists
and we have that
$$ \vert \int_0^T u(s) dW_s \vert \le G \Vert u \Vert_{\al,1},$$
where $G$ is a random variable with moments of all orders (see Lemma
7.5 in \cite{NR}). Moreover, if the trajectories of $u$ belong to
$W_0^{\al,\infty} (0,T)$, then the indefinite integral $\int_0^T u(s)
dW_s $ is H\"older continuous of order $1-\al$ and with trajectories
in $W_0^{\al,\infty} (0,T)$. As a simple consequence of these facts,
we get the following two proofs:

 \bigskip

\begin{dem2} {\text\bf Proof of Theorem \ref{teoexis}:} The existence and uniqueness of solution follows from
Theorem \ref{teoexi}. The existence of moment of any order is a consequence of
Lemma \ref{supnorminfprim}. Note only that if $\al < (2-\gamma)/4$ then $1/(1-\varphi(\gamma,\al)) < 2$ and
$ E( \exp(C \Lambda_\al (W)^{1/(1-\varphi(\gamma,\al))}
)) < \infty.$

\hfill $\Box$
\end{dem2}

\begin{dem2} {\text\bf Proof of Theorem \ref{convstoch}:} It suffices to apply Theorem
\ref{teoconver} to obtain the almost-sure convergence. The convergence in $L^p$ is obtained by a dominated
convergence argument since by Lemma \ref{supnorminf} we have that for any $r \in (-r_0,0)$
\begin{eqnarray*}
&&\Vert X- X^r \Vert_{\al,\infty} \le \Vert X \Vert_{\al,\infty}+ \Vert X^r \Vert_{\al,\infty(r)}
 \\&& \qquad  \le 2 C_\al^{(0)} (\Vert \eta
\Vert_{\al,\infty(-r_0,0)}+  \Lambda_\al(W) + 1) e^{C_\al^{(1)} +
C_\al^{(2)} \Lambda_\al(W)^{\frac{1}{1-\varphi(\gamma,\al)}}}:=Y,
\end{eqnarray*}
and $E(Y^p)< \infty$ for all $p \ge 1$.
\hfill $\Box$
\end{dem2}

\section{Appendix}

We recall two results from \cite{NR}: a fixed point theorem (see Lemma
7.2. page 75) and some algebraic computations, whose proof can be
easily derived from Lemma 7.1.

\begin{lema}\label{puntfixe}
Let $(X,\rho)$ be a complete metric space and $\rho_0, \rho_1, \rho_2$  some metrics
on $X$ equivalent to $\rho$.
Assume that $\cL: X \to X$ satisfies
\begin{enumerate}
\item there exists $\mu_0 > 0, x_0 \in X$ such that if $B_0 = \{x \in X:
\rho_0 (x_0, x) \le \mu_0 \}$
then  $\cL (B_0) \subset B_0$,
\item there exist $\varphi: (X, \rho) \to [0,+\infty]$ lower semicontinous function and some positive constants
$C_0, K_0$ such that denoting $N_\varphi (a) = \{ x \in X: \varphi(x) \le a \}$
\begin{enumerate}
\item $\cL (B_0) \subset N_\varphi(C_0)$,
\item $\rho_1 (\cL(x),\cL(y)) \le K_0 \rho_1 (x,y), \forall x,y \in N_\varphi(C_0) \cap B_0,$
\end{enumerate}
\item there exists $a \in (0,1)$ such that $\rho_2 ( \cL(x),\cL(y) ) \le a \rho_2 (x,y),
\forall x,y \in \cL(B_0)$.
\end{enumerate}
Then, there exists $x^* \in \cL (B_0) \subset X $ such that $x^*=\cL (x^*)$.
\end{lema}

\begin{lema}\label{calcDavid} Assume that $\si$ satisfies hypothesis {\bf (H1)}. Then for any
$f,g: [0,T] \to \R$ with $\Vert f \Vert_\infty \le N$ and $\Vert g \Vert_\infty \le N$  we have that
\begin{eqnarray*}
& &\int_0^t \frac{\vert \si(t,f(t))- \si(s,f(s)) - \si(t,h(t))+ \si(s,h(s)) \vert}{(t-s)^{\al +1}} ds \cr
& & \le M_0 \int_0^t \frac{\vert f(t)-f(s) - h(t)+ h(s) \vert}{(t-s)^{\al +1}} ds + \frac{M_0}{\beta-\al}
\vert f(t) - h(t) \vert t^{\beta - \al}\cr
& & \quad + M_N \vert f(t) - h(t) \vert
\left(\int_0^t \frac{\vert f(t)-f(s) \vert^\delta}{(t-s)^{\al +1}} ds +
\int_0^t \frac{\vert h(t)-h(s) \vert^\delta}{(t-s)^{\al +1}} ds \right).
\end{eqnarray*}

\end{lema}

\section*{Acknowledgements}
This work was partially supported by DGES Grants MTM2006-01351 (Carles Rovira).


\begin{thebibliography}{99}

\bibitem{FR}
M. Ferrante and C. Rovira: \it Stochastic delay
differential equations driven by fractional Brownian motion with
Hurst parameter $H > 1/2$. \rm Bernoulli {\bf 12} (2006), 85-100.


\bibitem{LT} J. Le\'on and S. Tindel: Private communication.



\bibitem{L} T. Lyons: \it Differential equations driven by rough signals (I): An extension
of an inequality of L. C. Young. \rm Mathematical Research Letters {\bf 1} (1994)
451-464 .

\bibitem{M} S.-E. A. Mohammed: \it Stochastic differential systems with memory: theory, examples and
applications. In Stochastic Analysis and Related Topics VI (L. Decreusefond, J. Gjerde, B. Øksendal
and A.S. \"Ust\"unel, eds), \rm Birkh\"auser, Boston, 1-77 (1998).


\bibitem{NNT}
A. Neuenkirch, I. Nourdin and S. Tindel:
\it Delay equations driven by rough paths.
\rm Electron. J. Probab.  {\bf 13 } (2008),  2031--2068.


\bibitem{NR} D. Nualart, A.  Rascanu: \it Differential equations driven by fractional Brownian
motion. \rm Collect. Math. {\bf 53} (2002) 55-81.


\bibitem{SKM} S.G. Samko, A.A. Kilbas and O. Marichev : \it Fractional Integrals and
Derivatives. Theory and Applications. \rm Gordon and Breach (1993) .


\bibitem{Y} L.C. Young: \it An inequality of the H\"older type connected with Stieltjes integration.
\rm Acta Math. {\bf 67} (1936) 251-282.

\bibitem{Z} M. Z\"ahle: \it Integration with respect to fractal functions and stochastic calculus.
I. \rm Prob. Theory Relat. Fields {\bf 111} (1998) 333-374.




\end{thebibliography}
\end{document}